\documentclass[a4paper,12pt]
{scrartcl}
\usepackage{typearea}
\typearea{13}

\usepackage{amsmath}
\usepackage{amssymb}
\usepackage[dvipdfmx]{graphicx}
\usepackage{float}
\usepackage{ulem}
\usepackage{enumitem}
\usepackage{ascmac}
\usepackage{framed}
\usepackage{comment}
\usepackage{here}
\usepackage[all]{xy}
\usepackage[b]{esvect}
\usepackage{color}
\usepackage{newtxtext,newtxmath}
\usepackage{multicol}

\usepackage{tabularx}
   \newcolumntype{C}{>{\centering\arraybackslash}X}
   \newcolumntype{L}{>{\raggedright\arraybackslash}X}
   \newcolumntype{R}{>{\raggedleft\arraybackslash}X}
   \newcolumntype{D}{>{\scriptsize\centering\arraybackslash}X}

\newtheorem{theorem}{Theorem}
\newtheorem{reftheorem}{Theorem}

\newtheorem{claim}{Claim}[section]

\newtheorem{refconjecture}[reftheorem]{Conjecture}
\newtheorem{remark}{Remark}

\newtheorem{Fact}{Fact}

\newcommand{\ora}{\vv}
\newcommand{\oraP}[1]{\vv{P}_{\!#1}\hspace{0.3mm}}

\newcommand{\sm}{\setminus}

\newcommand{\fl}{\flat}
\newcommand{\sh}{\sharp}
\newcommand{\el}{\ell}

\newcommand{\lp}{\el^{\sh}}

\newcommand{\proof}{\medbreak\noindent\textit{Proof.}\quad}
\newcommand{\qed}{{$\quad\square$\vspace{3.6mm}}}
\newcommand{\dist}{\textup{dist}}

\newcommand{\CF}{C}
\newcommand{\HF}{H}
\newcommand{\TF}{T}
\newcommand{\CT}{V_{c}}
\newcommand{\NCT}{V_{nc}}

\newcommand{\yb}[1]{\colorbox{yellow}{#1}}

\numberwithin{equation}{section}

\begin{document}

\title{Minimum degree conditions\\
for the existence of a sequence of cycles\\
whose lengths differ by one or two}

\author{Shuya Chiba\thanks{Applied Mathematics, Faculty of Advanced Science and Technology, 
Kumamoto University, 
2-39-1 Kurokami, Kumamoto 860-8555, Japan; 
E-mail address: \texttt{schiba@kumamoto-u.ac.jp}; 
This work was supported by JSPS KAKENHI Grant Number 20K03720}
\and 
Katsuhiro Ota\thanks{Department of Mathematics, Keio University,
3-14-1 Hiyoshi, Kohoku-ku, Yokohama 223-8522, Japan;
E-mail address: \texttt{ohta@math.keio.ac.jp};
This work was supported by JSPS KAKENHI Grant Number 16H03952}
\and
Tomoki Yamashita\thanks{Department of Science, Kindai University, 
3-4-1 Kowakae, Higashi-Osaka, Osaka 577-8502, Japan; 
E-mail address: \texttt{yamashita@math.kindai.ac.jp}; 
This work was supported by JSPS KAKENHI Grant Number 16K05262}
}

\date{}

\maketitle

\vspace{-24pt}
\begin{abstract}
Gao, Huo, Liu and Ma (2019) 
proved a result on the existence of paths connecting specified two vertices whose lengths differ by one or two.
By using this result, they settled
two famous conjectures due to Thomassen (1983). 
In this paper, 
we improve their result, 
and obtain 
a generalization of a result of Bondy and Vince (1998).

\medskip
\noindent
\textit{Keywords}: Cycle length, Path length, Minimum degree
\\
\noindent
\textit{AMS Subject Classification}: 05C38
\end{abstract}

\section{Introduction}
\label{sec:introduction}

All graphs considered in this paper are finite undirected graphs without loops or multiple edges. 
For a graph $G$, 
we denote by 
$V(G)$ and $E(G)$ the vertex set and the edge set of $G$, respectively, 
and 
$\deg_{G}(v)$ denotes the degree of a vertex $v$ in $G$.

In 1983, 
Thomassen proposed the following two conjectures.

\begin{refconjecture}[Thomassen \cite{T1983}]
\label{T1983 all even length k}
For a positive integer $k$, 
every graph of minimum degree at least $k + 1$ 
contains cycles of all even lengths modulo $k$. 
\end{refconjecture}

\begin{refconjecture}[Thomassen~\cite{T1983}]
\label{T1983 all length k}
For a positive integer $k$, 
every $2$-connected non-bipartite graph 
of minimum degree at least $k + 1$ 
contains cycles of all lengths modulo $k$. 
\end{refconjecture}

The above conjectures originated from the conjecture of Burr and Erd\H{o}s 
concerning the extremal problem for the existence of cycles with prescribed lengths modulo $k$ 
(see \cite{Erdos1976}). 
We refer the reader to \cite{LM2018} for more details. 
In 2018, 
Liu and Ma proved that Conjectures~\ref{T1983 all even length k} and \ref{T1983 all length k} are true for all even integers~$k$ 
by considering the existence of a sequence of paths whose lengths differ by two 
in bipartite graphs, 
see \cite[Theorem~1.9]{LM2018}.

Recently, 
Gao, Huo, Liu and Ma~\cite{GHLM2019+} 
announced that 
they had confirmed Conjectures~\ref{T1983 all even length k} and \ref{T1983 all length k} for 
all integers $k$ 
by using the following theorem. 
Here, 
we say that 
a sequence of $k$ paths (or $k$ cycles) $P_{1}, \dots, P_{k}$ is 
\textit{admissible}
if $|E(P_{1})| \ge 2$, and 
either $|E(P_{i+1})| - |E(P_{i})| = 1$ 
for $1 \le i \le k - 1$ 
or
$|E(P_{i+1})| - |E(P_{i})| = 2$ 
for $1 \le i \le k - 1$.

\begin{reftheorem}[Gao et al.~{\cite[Theorem 1.2]{GHLM2019+}}]
\label{thm:Gao et al}
Let $k$ be a positive integer, 
and let 
$G$ be a $2$-connected graph,
and $x,y$ be two distinct vertices of $G$. 
If 
$\deg_{G}(v) \ge k+1$
for each $v \in V(G) \setminus \{x, y\}$, 
then 
$G$ contains $k$ admissible paths from $x$ to $y$. 
\end{reftheorem}

In this paper, 
we show that 
the degree condition in Theorem~\ref{thm:Gao et al} can be relaxed as follows.

\begin{theorem}\label{thm:three vertices} 
Let $k$ be a positive integer, 
and let 
$G$ be a $2$-connected graph, 
$x,y$ be two distinct vertices of $G$, 
and $z$ be a vertex of $G$ (possibly $z \in \{x, y\}$) 
such that $V(G) \setminus \{x, y, z\} \neq \emptyset$. 
If 
$\deg_{G}(v) \ge k+1$ for each $v \in V(G) \setminus \{x, y, z\}$, 
then 
$G$ contains $k$ admissible paths from $x$ to $y$. 
\end{theorem}

This study also originated from the question of 
whether every graph of minimum degree at least three contains two admissible cycles, 
which was raised by Erd\H{o}s (see \cite{BV1998}). 
In 1998, Bondy and Vince answered this queston by proving the following stronger theorem.

\begin{reftheorem}[Bondy, Vince~\cite{BV1998}]
\label{thm:BV1998}
Every graph of order at least three, 
having at most two vertices of degree less than three, contains two admissible cycles. 
\end{reftheorem}

They also conjectured that 
every graph of sufficiently large order, 
having at most $m$ vertices of degree less than three, 
contains two admissible cycles, and they gave some remarks for the case of small~$m$. 
In 2020, Gao and Ma~\cite{GM2020} settled the conjecture in the affirmative for all $m$.

We give the following another generalization of Theorem~\ref{thm:BV1998} 
by using Theorem~\ref{thm:three vertices}.

\begin{theorem}\label{thm:generalization of BV} 
For an integer $k \ge 2$, 
every graph of order at least three, 
having at most two vertices of degree less than $k + 1$, contains $k$ admissible cycles.
\end{theorem}

Note that 
a weaker version of Theorem~\ref{thm:generalization of BV} is obtained from 
Theorem~\ref{thm:Gao et al} (see~\cite[Theorem 1.3]{GHLM2019+}), 
and the result (and also Theorem~\ref{thm:generalization of BV}) 
settles Conjectures~\ref{T1983 all even length k} and \ref{T1983 all length k} for 
all integers~$k$.

\medskip
To show Theorem~\ref{thm:three vertices}, 
in the next section, 
we consider the existence of admissible paths in ``rooted graphs'' 
and give a stronger result than Theorem~\ref{thm:three vertices} (see Theorem~\ref{thm:three vertices in rooted graph} in Section~\ref{sec:preliminaries}). 
We also extend the concept of ``cores'' which were used in the argument of \cite{LM2018, GHLM2019+} 
in preparation for the proof of Theorem~\ref{thm:three vertices in rooted graph}. 
In Section~\ref{sec:proof of main}, we prove Theorem~\ref{thm:three vertices in rooted graph}  
and also give the proof of Theorem~\ref{thm:generalization of BV} 
at the end of Section~\ref{sec:proof of main}.

\section{Preliminaries}
\label{sec:preliminaries}

\subsection{Admissible paths in rooted graphs}

Let $G$ be a graph. 
A \textit{cut-vertex} of $G$ is a vertex whose removal increases the number of components of $G$. 
A \textit{block} of $G$ is a maximal connected subgraph of $G$ which has no cut-vertex, 
and 
a block $B$ of $G$ is called an \textit{end-block} if $B$ has at most one cut-vertex of $G$. 
If $G$ itself is connected and has no cut-vertex, then $G$ is a block and is also an end-block. 
For distinct vertices $x$ and $y$ of $G$, 
$(G, x, y)$ is called a \textit{rooted graph}. 
A rooted graph $(G, x, y)$ is \textit{$2$-connected} 
if 
\begin{enumerate}[label=(\texttt{R}\arabic*)]
\item
\label{BT=path} 
$G$ is a connected graph of order at least three with at most two end-blocks, 
and 
\item 
\label{xyEnd}
every end-block of $G$ contains at least one of $x$ and $y$ 
as a non-cut-vertex of $G$. 
\end{enumerate}
Note that $(G, x, y)$ is $2$-connected if and only if $G + xy$ 
(i.e., the graph obtained from $G$ 
by adding the edge $xy$ if $xy \notin E(G)$) is $2$-connected. 
We denote by $(G, x, y; z)$
a rooted graph $(G, x, y)$ with a specified vertex $z$ 
(this includes the case where $z \in \{x, y\}$ or $z \not\in V(G)$). 
For a rooted graph $(G, x, y; z)$, 
we define $\delta(G, x, y; z) = \min \{\deg_{G}(v) : v \in V(G) \setminus \{x, y,z \}\}$ 
if $V(G) \setminus \{x, y, z\} \neq \emptyset$; 
otherwise, let $\delta(G, x, y; z) = - \infty$.

\bigskip

In this paper,
instead of proving Theorem \ref{thm:three vertices},
we prove the following stronger theorem.

\begin{theorem}\label{thm:three vertices in rooted graph} 
Let $k$ be a positive integer, 
and let 
$(G, x, y; z)$ be a $2$-connected rooted graph. 
If $\delta(G, x, y; z) \ge k+1$,
then $G$ contains $k$ admissible paths from $x$ to $y$.
\end{theorem}

\subsection{Terminology and notation}
\label{subsec:terminology and notation}

In this subsection,  
we prepare terminology and notation which will be used in the proof of Theorem~\ref{thm:three vertices in rooted graph}.

Let $G$ be a graph. 
We denote by $N_{G}(v)$ the neighborhood of a vertex $v$ in $G$. 
For $S \subseteq V(G)$, 
we define  
$N_{G}(S) = \big( \bigcup_{v \in S}N_{G}(v) \big) \setminus S$. 
For $S \subseteq V(G)$, 
$G[S]$ denotes the subgraph of $G$ induced by $S$, 
and let $G - S = G[V(G) \setminus S]$. 
We denote by $\dist_{G}(u, v)$ the length of a shortest path from a vertex $u$ to a vertex $v$ in $G$. 
For $U, V \subseteq V(G)$ with $U \cap V = \emptyset$, 
a path in $G$ is a \textit{$(U, V)$-path} 
if 
one end-vertex of the path belongs to $U$, 
the other end-vertex belongs to $V$, 
and the internal vertices do not belong to $U \cup V$. 
We write a path $P$ with a given orientation as $\ora{P}$. 
For an oriented path $\ora{P}$ and $u, v \in V(P)$, 
a path from $u$ to $v$ along $\ora{P}$ is denoted by $u\ora{P}v$. 
For $t \ (\ge 2)$ pairwise vertex-disjoint sets $V_{1},\ldots, V_{t}$ of vertices, 
we define the graph $V_{1} \vee \cdots \vee V_{t}$ from the union of $V_{1}, \ldots, V_{t}$ 
by joining every vertex of $V_{i}$ to every vertex of $V_{i+1}$ for $1\leq i \leq t-1$. 
For convenience, 
we let
$V_{1} \vee \cdots \vee V_{t} \vee \emptyset = V_{1} \vee \cdots \vee V_{t}$.

\begin{figure}[H]
\begin{center}
{\includegraphics[pagebox=artbox]{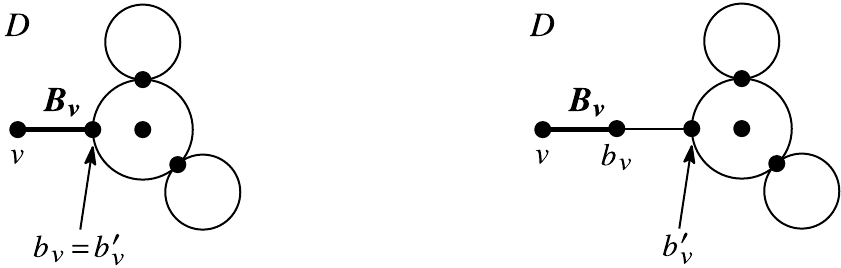}}
\end{center}
\caption{The $v$-end-block $B_{v}$ of $D$}
\label{vendblockFIG}
\end{figure}

Let 
$D$ be a connected graph 
and 
$v$ be a vertex. 
The \textit{$v$-end-block of $D$} 
is an end-block $B_{v}$ with cut-vertex $b_{v}$ in $D$ 
such that 
$V(B_{v}) = \{v, b_{v}\}$. 
The $v$-end-block of $D$, if exists, is unique, 
and so we always denote it by $B_{v}$ for a vertex $v$. 
We also denote by $b_{v}$ the unique cut-vertex of $D$ which is contained in $B_{v}$. 
If $\deg_{D}(b_{v}) = 2$, 
then let $b_{v}'$ denote the unique neighbor of $b_{v}$ in $D$ which is not $v$; 
otherwise, let $b_{v}' = b_{v}$. See Figure~\ref{vendblockFIG}.

Throughout the rest of this paper, 
we often denote the singleton set $\{v\}$ by $v$, 
and 
we often identify a subgraph $H$ of $G$ with its vertex set $V(H)$.

\subsection{The concept of cores}
\label{subsec:core}

In this subsection, 
we extend the concept of cores which were used in the argument of \cite{LM2018, GHLM2019+}.

Let $\el$ be an integer. Let $x$ be a vertex of a graph $G$, and let $H$ be a subgraph of $G$. 
\begin{itemize} 
\item 
$H$ is called an \textit{$\el$-core of type~1 with respect to $x$} 
if 
$H = x \vee T \vee S$,
where 
$\el \ge 1$, 
$S = \emptyset$ 
and 
$T$ is a clique of order exactly $\el+1$ in $G$. 
\item 
$H$ is called an \textit{$\el$-core of type~2 with respect to $x$} 
if 
$H = x \vee S \vee T$,
where 
$\el \ge 2$, 
$S$ is an independent set of order exactly $2$
and
$T$ is a clique of order exactly $\el$.
\item 
$H$ is called an \textit{$\el$-core of type~3 with respect to $x$}
if 
$H = x \vee T \vee S$,
where 
$\el \ge 0$ and, 
$S$ and $T$ 
are independent sets of orders exactly $\el$ and at least $\max\{\el+1,2\}$,
respectively. 
(Since $\el \ge 0$, $S$ may be an empty set.) 
\end{itemize}

\begin{figure}[H]
\begin{center}
{\includegraphics[pagebox=artbox,scale=1]{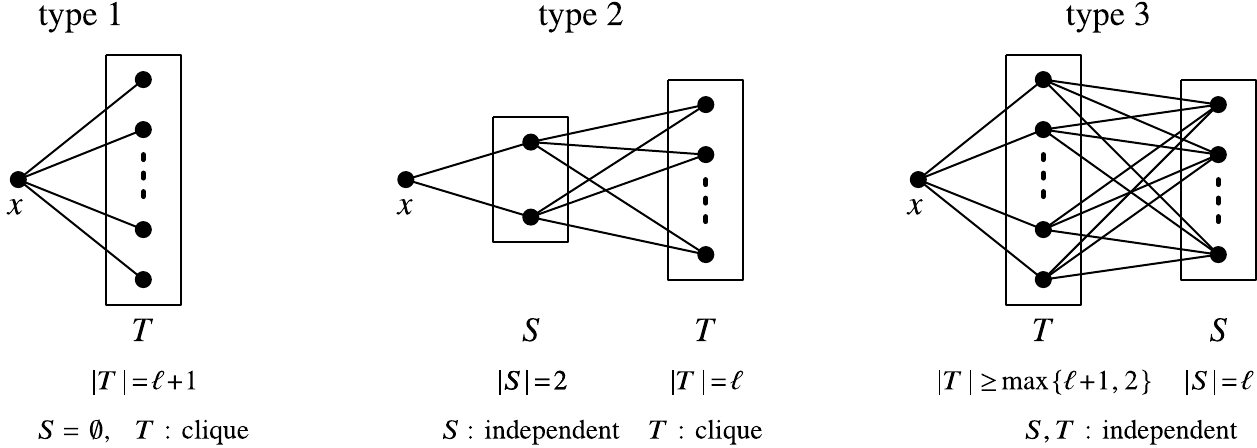}}
\end{center}
\caption{The structures of $\el$-cores of types~1, 2 and 3 with respect to $x$}
\label{coreFIG}
\end{figure}

\noindent
See Figure~\ref{coreFIG}. We say that $H$ is an $\el$-core with respect to $x$ 
when there is no need to specify the type. 
We also say that 
$H$ is an \textit{$\el$-core with respect to $(x, y)$} 
if $H$ is an $\el$-core with respect to $x$, and $y$ is a vertex of $V(G) \setminus V(H)$. 
In what follows, ``a core'' always means an $\el$-core for some integer $\el$.

\begin{remark}
\label{remark:lcore}
Let $G$ be a graph, and $x, y$ be two distinct vertices of $G$. 
If $\deg_{G}(x) \ge 2$ and $xy \notin E(G)$, 
then 
there exists a core of type~1 or type~3 with respect to $(x, y)$ in $G$. 
\end{remark}

We give three facts which will be used frequently
in the proof of Theorem~\ref{thm:three vertices in rooted graph}. 
Here, 
an admissible sequence 
which changed the condition $|E(P_{1})| \ge 2$ into $|E(P_{1})| \ge 1$, 
is said to be 
\textit{semi-admissible}.

\begin{Fact}[Gao et al.~{\cite[Lemma 3.2]{GHLM2019+}}]
\label{fact:Lemma3.2} 
Let $s, t$ be positive integers. 
Let $G$ be a graph, $x, y$ be two distinct vertices 
and $U \subseteq V(G) \setminus \{x, y\}$. 
Assume that $G$ contains $s$ semi-admissible $(U, y)$-paths $P_{1}, \dots, P_{s}$, 
and let $u_{i}$ be the unique vertex of $V(P_{i}) \cap U$ for $1 \le i \le s$. 
Assume further that 
for each $1 \le i \le s$, 
$G - (V(P_{i}) \setminus \{u_{i}\})$ contains 
$t$ semi-admissible $(x, u_{i})$-paths $Q_{i, 1}, \dots, Q_{i, t}$. 
If $|V(Q_{1,j})| = |V(Q_{2,j})| = \cdots = |V(Q_{s, j})|$ for $1 \le j \le t$, 
then 
$G$ contains $s + t - 1$ admissible $(x, y)$-paths.
\end{Fact}

\begin{Fact}
\label{fact:paths in H} 
Let $H$ be an $\el$-core with respect to $x$. 
Then the following hold. 
\begin{enumerate}[label=(\arabic*)]
\item 
\label{xS}
If $H$ is of type~2, 
then for any $s \in S$,
$H$ contains $\el$ admissible $(x, s)$-paths of lengths $3, 4, \ldots, \el+2$; 
if $H$ is of type~3, 
then for any $s \in S$ and $t \in T$, 
$H-t$ contains $\el$ admissible $(x, s)$-paths of lengths $2,4,\ldots,2\el$.
\item 
\label{xT}
For any $t \in T$,
$H$ contains $\el+1$ semi-admissible $(x, t)$-paths
of lengths $1,2,\ldots,\el+1$ (if $H$ is of type 1) 
or $2,3,\ldots,\el+2$ (if $H$ is of type 2) 
or $1,3,\ldots,2\el+1$ (if $H$ is of type 3). 
In particular, if $H$ is of type~3 and $|T| \ge \el + 2$, 
then for any $t, t' \in T$ with $t \neq t'$, 
$H-t'$ contains $\el + 1$ semi-admissible $(x, t)$-paths of lengths $1,3,\ldots,2\el+1$. 
\end{enumerate}
\end{Fact}

Fact~\ref{fact:paths in H} immediately yields the following.

\begin{Fact}
\label{fact:l >= k or l >= k-1} 
Let $k$ be a positive integer 
and 
$(G, x, y)$ be a $2$-connected rooted graph. 
Let 
$H$ be an $\el$-core with respect to $(x, y)$ 
and $C$ be the component of $G-V(H)$ such that $y \in V(C)$. 
Assume that $G$ does not contain $k$ admissible $(x, y)$-paths. 
Then 
(1) $\el \le k -1$, 
and 
(2) if 
$N_{G}(C) \cap T \neq \emptyset$, 
then $\el \le k-2$.
\end{Fact}

\section{Proof of Theorem~\ref{thm:three vertices in rooted graph}}
\label{sec:proof of main}

\noindent
\textit{Proof of Theorem~\ref{thm:three vertices in rooted graph}}.
We prove it by induction on $|V(G)| + |E(G)|$.
Let $(G, x, y; z)$ be a minimum counterexample with respect to $|V(G)| + |E(G)|$.

\begin{claim}
\label{claim:mindeg3andorder5} 
(1) $k \ge 2$ (and so $\delta(G, x, y; z) \ge 3$), 
and 
(2) $|V(G)| \ge 5$. 
\end{claim}
\proof 
(1) If $k = 1$, then by \ref{BT=path} and \ref{xyEnd}, 
we can easily see that $G$ contains an $(x, y)$-path of length at least $2$, a contradiction. 
Thus $k \ge 2$, 
and so $\delta(G, x, y; z) \ge k + 1 \ge 3$. 

(2) Since $\delta(G, x, y; z) \ge k + 1 \ge 3$, 
we have $|V(G)| \ge 4$. 
Suppose that $|V(G)| = 4$. 
Note that, then $\delta(G, x, y; z) = k + 1 = 3$. 
Let $u$ and $v$ be two distinct vertices of $V(G) \setminus \{x, y\}$ 
such that $u \neq z$. 
Since $\deg_{G}(u) = 3$, 
we have $N_{G}(u) = \{x, y, v\}$. 
By \ref{xyEnd}, 
we also have $N_{G}(v) \cap \{x, y\} \neq \emptyset$, 
say $xv \in E(G)$ up to symmetry, 
and then 
$xuy$ and $xvuy$ are $k \ (=2)$ admissible $(x, y)$-paths, 
a contradiction.
\qed

\begin{claim}
\label{claim:2-connected and xy notin E} 
(1) $G$ is $2$-connected 
and 
(2) $\{x,y,z\}$ is independent.
\end{claim}
\proof 
(1) Suppose that $G$ is not $2$-connected. 
Then by \ref{BT=path}, 
$G$ has a cut vertex $c$ 
and $G-c$ has exactly two components $C_{1}$ and $C_{2}$. 
By \ref{xyEnd}, without loss of generality,
we may assume that $x \in V(C_{1})$ and $y \in V(C_{2})$.
Since $V(G)\sm\{x,y,z,c\}\not=\emptyset$ by Claim~\ref{claim:mindeg3andorder5}~(2), 
and by the symmetry of $x$ and $y$, 
we may assume that $V(C_{1}) \setminus \{x, z\} \neq \emptyset$.
Let $G_{i} = G[C_{i} \cup c]$ for $i \in \{1, 2\}$. 
Then
$(G_{1}, x, c; z)$ is a 2-connected rooted graph 
such that 
$\delta(G_{1}, x, c; z) \ge \delta(G, x, y; z)$.
Hence by the induction hypothesis, 
$G_{1}$ contains $k$ admissible $(x, c)$-paths $\oraP{1}, \dots, \oraP{k}$.
Let $\ora{Q}$ be a $(c, y)$-path in $G_{2}$.
Then 
$x\oraP{i} c \ora{Q} y$ ($1 \le i \le k$)
are $k$ admissible $(x, y)$-paths in $G$, 
a contradiction. 

(2) 
Suppose that $xv \in E(G)$ 
for some $v \in \{y,z\}$, 
and choose such a vertex $v$ so that $v = y$ if possible. 
If 
$G - xv$ 
(i.e., the graph obtained from $G$ by deleting the edge $xv$) 
is $2$-connected, 
then 
by the induction hypothesis, 
it follows that 
$G - xv$ (and also $G$) contains $k$ admissible $(x, y)$-paths, 
a contradiction. 
Thus 
$G - xv$ is not $2$-connected. 
Since 
$G$ is $2$-connected by Claim~\ref{claim:2-connected and xy notin E}~(1), 
this implies that 
$(G-xv, x, v)$ is a $2$-connected rooted graph with exactly two end-blocks.

Let 
$B_{1}, \dots, B_{t}$ ($t \ge 2$) be all the blocks of $G -xv$ 
such that 
$V(B_{i}) \cap V(B_{i + 1}) \neq \emptyset$ for $1 \le i \le t-1$, 
say $V(B_{i}) \cap V(B_{i + 1}) = \{b_{i}\}$ for $1 \le i \le t-1$. 
Without loss of generality, 
we may assume that 
$x \in V(B_{1}) \setminus \{b_{1}\}$ 
and 
$v \in V(B_{t}) \setminus \{b_{t-1}\}$. 
Then 
$y \in V(B_{p}) \setminus \{b_{p-1}\}$
for some $p$ with $1 \le p \le t$, 
where we let $b_{0} = x$.

Suppose that $p=t$.
Then
$(G - xv, x, y; z)$ is a $2$-connected rooted graph 
such that 
$\delta(G - xv, x, y; z) = \delta(G, x, y; z)$.  
Hence, by the induction hypothesis,
$G -xv$ (and also $G$) contains $k$ admissible  $(x, y)$-paths, 
a contradiction. 
Thus $p \le t - 1$. 
This implies that $v \neq y$, that is, $v = z$. 
Then by the choice of $v$, we have $xy \notin E(G)$.

Let 
$G' = G[\bigcup_{1 \le i \le p}V(B_{i})]$,
and let $z' = b_{p}$. 
Note that $z \notin V(G')$. 
Note also that 
if $p = 1$, 
then since $xy \not\in E(G)$, 
$V(G') \setminus \{x, y, z'\} = V(B_{1}) \setminus \{x, y, b_{1}\} \neq \emptyset$ holds; 
if $p \ge 2$, 
$V(G') \setminus \{x, y, z'\} \neq \emptyset$ clearly holds. 
Then 
$(G', x, y; z')$ is a $2$-connected rooted graph 
such that $\delta(G', x, y; z') \ge \delta(G, x, y; z)$, 
and so the 
the induction hypothesis yields that 
$G'$ (and also $G$) contains $k$ admissible $(x, y)$-paths, 
a contradiction. 
Thus $xv \notin E(G)$ for each $v \in \{y, z\}$.
By the symmetry of $x$ and $y$, we also have $yz \notin E(G)$. 
\qed

By Remark~\ref{remark:lcore} 
and 
Claim~\ref{claim:2-connected and xy notin E},
there exist 
cores with respect to $(x, y)$ and $(y, x)$, respectively, in $G$. 
By the symmetry of $x$ and $y$, 
we can rename the vertices $x$ and $y$ so that 
\begin{enumerate}[label=(\texttt{X\!Y}\arabic*)]
\item
\label{refno}
there exists a core $H$ with respect to $(x, y)$ so that 
the number of type of $H$ 
 is as small as possible, 
\item
\label{degxy}
$\deg_{G}(x) \le \deg_{G}(y)$, subject to \ref{refno},  and 
\item
\label{dist}
$\dist_{G}(x, z) \le \dist_{G}(y, z)$, subject to \ref{refno} and \ref{degxy}. 
\end{enumerate}

\medskip
Let $H$ be an $\el$-core with respect to $(x, y)$ in $G$ for some integer $\el$, 
and let $C$ be the component of $G- V(H)$ such that $y \in V(C)$. 
Choose $H$ so that 
\begin{enumerate}[label=(\texttt{H}\arabic*)]
\item
\label{renfomin}
the number of type of $H$ 
is as small as possible, and 
\item
\label{Hlarge} 
subject to \ref{renfomin}, 
\begin{enumerate}[label=(\texttt{H}2-\arabic*)]
\item
\label{Tlarge} 
if $H$ is of type 1 or type 2, 
then $|T|$ is maximum;
\item 
\label{TSlarge} 
if $H$ is of type 3, 
then 
(i) $|S|$ is maximum,
(ii) $|T|$ is maximum, subject to (i). 
\item
\label{TSmodified} 
If $H$ and $C$ satify the following condition~\ref{transform}:
\begin{enumerate}[label=(\texttt{T})]
\item
\label{transform}
$H$ is of type~3, 
$|T| \ge 3$, 
$V(C)=\{y\}$, 
$N_{G}(x)=N_{G}(y)=T$,
and there exists a component $D_{0}$ 
of $G - V(H)$ 
such that $D_{0}\not=C$, $V(D_{0}) \setminus \{z\} \neq \emptyset$
and $N_{G}(D_{0}) \cap T\not=\emptyset$,
\end{enumerate}
then
let 
$t_{0} \in N_{G}(D_{0}) \cap T \ ({}\cap N_{G}(y))$, 
and 
we modify $H$ (and $C$ depending on $H$) 
by resetting $\el$, $S$ and $T$ as follows:
\begin{enumerate}[label=(\texttt{M}\arabic*)]
\item 
\label{M1}
if $|T| = |S| + 1$\footnote{Note that, in this case, $\el \ge 2$.}, then
let $s_{0} \in S$, 
and we reset 
$\el := \el - 1$, 
$S :=S \sm \{s_{0}\}$
and 
$T :=T \sm \{t_{0}\}$; 
\item 
\label{M2}
if $|T| \ge |S|+2$, then 
we reset 
$\el := \el$, 
$S := S$ 
and 
$T:=T \sm \{t_{0}\}$.
\end{enumerate}
\end{enumerate}
\end{enumerate}
\begin{figure}[H]
\begin{center}
{\includegraphics[pagebox=artbox]{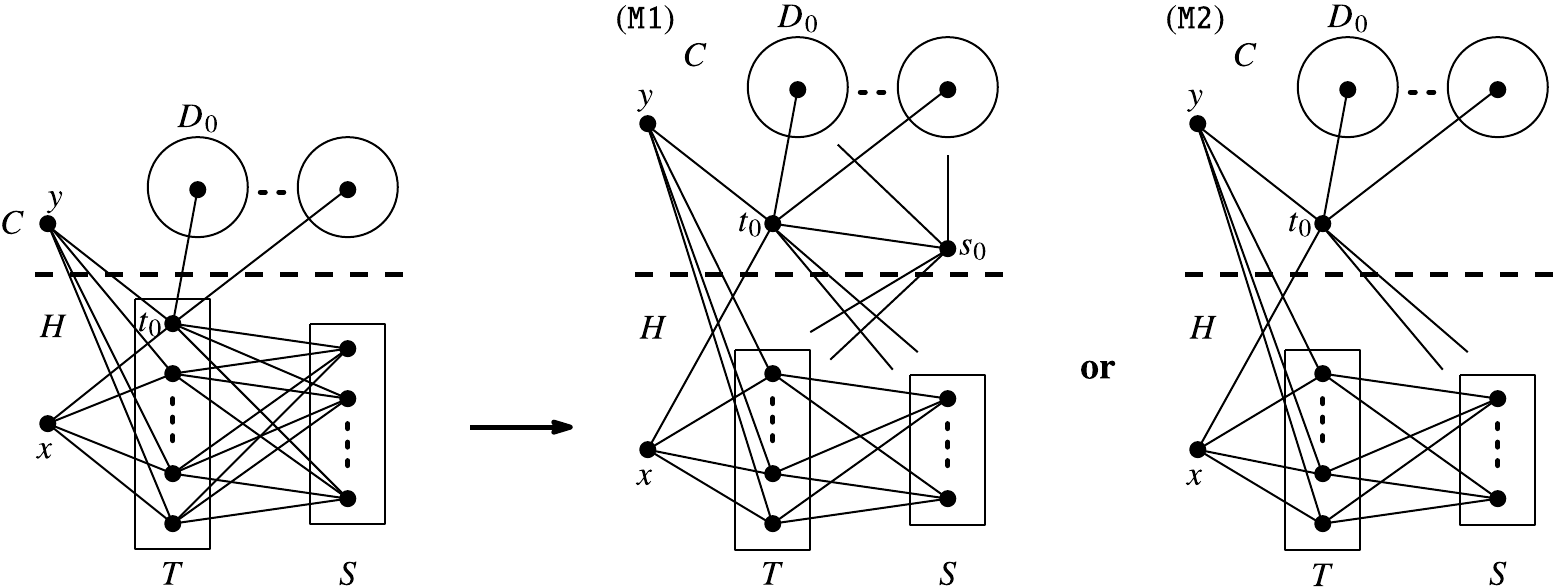}}
\end{center}
\caption{The modifications of $H$ and $C$}
\label{transformFIG}
\end{figure}
Note that 
if we modify $H$ as in \ref{TSmodified}, 
then 
the new graph $H$ is also an $\el$-core of type 3 with respect to $(x, y)$ 
(see also Figure~\ref{transformFIG}). 
However, to make the difference clear,
we sometimes say that
$H$ is of type~3$^{\fl}$ if $H$ is modified as above; 
otherwise 
$H$ is of type~3$^{\natural}$.

\begin{claim}
\label{claim:at most l+1}
If $v$ is a vertex of $V(G) \sm (V(H) \cup \{y,t_{0}\})$, 
then 
$|N_{G}(v) \cap V(H)| \le \el+1$, 
where $t_{0} = y$ for the case where $H$ is not of type~3$^{\fl}$.
In particular, 
if the equality holds,
then
$N_{G}(v) \cap T \not=\emptyset$.
%
\end{claim}
\proof 
Let $v$ be a vertex of $V(G) \sm (V(H) \cup \{y, t_{0}\})$. We show the claim as follows.

Assume first that $H$ is of type~1. 
If $|N_{G}(v) \cap V(H)| \ge \el + 2$, 
then 
we have $N_{G}(v) \cap V(H) = x \cup T$, 
which contradicts the maximality of $T$ (see \ref{Tlarge}). 
Thus $|N_{G}(v) \cap V(H)| \le \el + 1$. 
If the equality holds, 
we clearly have 
$N_{G}(v) \cap T \not=\emptyset$, since $\el + 1 \ge 2$.

Assume next that $H$ is of type~2, 
and suppose that 
$|N_{G}(v) \cap V(H)| \ge \el + 2$.
Since there exists no core of type~1 with respect to $(x, y)$ by \ref{renfomin}, 
we have $x \not\in N_{G}(v)$ or $N_{G}(v) \cap S =\emptyset$. 
Since $|N_{G}(v) \cap V(H)| \ge \el + 2$, 
$|S| = 2$ and $|T| = \el$, 
this yields that 
$N_{G}(v) \cap V(H) = S \cup T$, 
which contradicts the maximality of $T$ (see \ref{Tlarge}). 
Thus $|N_{G}(v) \cap V(H)| \le \el + 1$. 
If the equality holds, 
then 
since 
$x \not\in N_{G}(v)$ or $N_{G}(v) \cap S =\emptyset$ holds, 
it follows that 
$|N_{G}(v) \cap T| = (\el + 1) - |N_{G}(v) \cap (x \cup S)| \ge (\el + 1) - 2 = \el -1 \ge 1$. 
Thus 
we have $N_{G}(v) \cap T \not=\emptyset$.

Assume finally that $H$ is of type~3. 
Suppose that 
$|N_{G}(v) \cap V(H)| \ge \el + 1$.
Since there exists no core of type~1 with respect to $(x, y)$ by \ref{renfomin},
we have $x \not\in N_{G}(v)$ or $N_{G}(v) \cap T =\emptyset$.
Since there also exists no core of type~2 with respect to $(x, y)$ by \ref{renfomin},
we have $|N_{G}(v)\cap T| \le 1$ or $N_{G}(v) \cap S =\emptyset$. 
If $H$ is of either type~3$^{\natural}$ or type~3$^{\fl}$ in \ref{M2}, 
then $|N_{G}(v) \cap T| \le \el+1$ (by \ref{TSlarge}-(i)); 
if $H$ is of type~3$^{\fl}$ in \ref{M1},
then we have $|N_{G}(v) \cap T| \le |T| = \el+1$.
If $H$ is of type~3$^{\natural}$,
then
$x \cup S \not\subseteq N_{G}(v)$ (by \ref{TSlarge}-(ii));
if $H$ is of type~3$^{\fl}$,
then
since $N_{G}(x)=T \cup \{t_{0}\}$ and $v\not=t_{0}$,
we have $x \cup S \not\subseteq N_{G}(v)$.
Since $|S|=\el$,
combining these facts 
yields that
$|N_{G}(v) \cap V(H)| = \el+1$,
and
$N_{G}(v) \cap T \neq \emptyset $. 
\qed

By Claim~\ref{claim:2-connected and xy notin E}~(2) 
and \ref{transform}, 
we can easily obtain the following.

\begin{claim}
\label{claim:transformed core}
If $H$ is of type 3$^{\fl}$,
then 
$|V(\CF) \setminus \{y, z\}| \ge 2$.
\end{claim}
\proof 
Since $yz \notin E(G)$ by Claim~\ref{claim:2-connected and xy notin E}~(2) 
and 
$V(D_{0}) \setminus \{z\} \neq \emptyset$
by \ref{transform}, 
we have 
$|V(\CF) \setminus \{y, z\}| \ge 
 |V(D_{0}) \setminus \{z\}| 
+ |\{t_{0}\}| \ge 2$.
\qed

We now divide the proof into two cases according to $V(\CF) = \{y\}$ or $V(\CF) \neq \{y\}$.

\bigskip
\noindent\textbf{Case~1.}
$V(\CF) = \{y\}$. 

Note that by Claim~\ref{claim:transformed core}, 
$H$ is not of type 3$^{\fl}$.

\begin{claim}
\label{claim:neighbor of x and y}
If $H$ is of type 1 or type 3,
then 
$N_{G}(x) = N_{G}(y)=T$.
If $H$ is of  type 2,
then 
$N_{G}(x) = N_{G}(y)=S$. 
\end{claim}
\proof 
Note that 
by Claim~\ref{claim:2-connected and xy notin E}, 
$|N_{G}(y) \cap V(H-x)| = \deg_{G}(y) \ge 2$. 

Assume first that 
$H$ is of type~1. 
Then 
$|N_{G}(y)\cap T| \ge 2$, 
and so 
there exists a core of type~1 
with respect to $(y, x)$. 
Since 
$N_{G}(y) \subseteq T \subseteq N_{G}(x)$, 
it follows from \ref{degxy} that 
$N_{G}(x) = N_{G}(y)= T$. 
Thus the claim follows.

Assume next that $H$ is of type~2. 
Since $|N_{G}(y) \cap V(H-x)| \ge 2$, 
and since 
there exists no core of type~1 with respect to $(y, x)$ by \ref{refno}, 
we have $N_{G}(y) = S$. 
This in particular implies that 
$y \vee S \vee T$ is an $\el$-core of type~2 with respect to $(y, x)$. 
Since $N_{G}(y) = S \subseteq N_{G}(x)$, 
it follows from \ref{degxy} that 
$N_{G}(x) = N_{G}(y) = S$. 
Thus the claim follows.

Assume finally that $H$ is of type~3. 
By \ref{refno} and \ref{renfomin}, 
there exist no cores of type 1 or type~2 with respect to $(y, x)$, 
and so 
any core with respect to $(y, x)$ is of type $3$ 
(by Remark~\ref{remark:lcore}, Claim~\ref{claim:2-connected and xy notin E}). 
This also implies that 
$N_{G}(y) \subseteq T$ or $N_{G}(y) \subseteq S$. 
Since 
$|T| \ge \max\{\el + 1, 2\} > \el = |S|$ and $T \subseteq N_{G}(x)$, 
it follows from \ref{degxy} that 
$N_{G}(x) = N_{G}(y)= T$. 
Thus the claim follows. 
\qed

\begin{claim}
\label{claim:notsubsetT}
Assume that 
$H$ is of either type~1 or type~3.
Let $D$ be a component of $G - V(H)$
such that $D \not=C$ and $V(D) \setminus \{z\} \neq \emptyset$.
Then 
$N_{G}(D) \cap S\not=\emptyset$. 
(This in particular implies that, 
if $H$ is of type~1, then 
$G - V(H)$ does not have a component $D$ 
such that $D \not=C$ and $V(D) \setminus \{z\} \neq \emptyset$.) 
\end{claim}
\proof
Suppose that 
$N_{G}(D) \subseteq T$. 
By Claim \ref{claim:neighbor of x and y},
there exists a vertex $t_{c\hspace*{-0.25mm}d} \in N_{G}(y) \cap N_{G}(D) \cap T$. 
Let $D^*$ be the graph obtained from $G[D \cup N_{G}(D)]$
by contracting $N_{G}(D) \sm\{t_{c\hspace*{-0.25mm}d}\}$
into a new vertex $t^{*}$. 
Since $N_{G}(D) \subseteq T$ and 
$G$ is $2$-connected, 
it follows that $(D^*, t^{*}, t_{c\hspace*{-0.25mm}d}; z)$ is a $2$-connected rooted graph. 
Since $V(D) \setminus \{z\} \neq \emptyset$, 
we also have $\emptyset \neq V(D^{*}) \setminus \{t^{*}, t_{c\hspace*{-0.25mm}d}, z\} \subseteq V(G) \setminus \{x, y, z\}$.
Let 
\[
\epsilon=
\begin{cases}
1 & \text{if $|T|=\el+1$,}\\
0 & \text{if $|T|\ge \el+2$.}
\end{cases}
\]
If 
$H$ is of type~3, 
then 
$|T| \ge \max \{\el + 1, 2\}$, and so $\el \ge 1$ holds for the case of $|T| = \el + 1$, 
which implies that $\el - \epsilon \ge 0$; 
if 
$H$ is of type~1, 
then since 
$|T| = \el + 1 \ge 2$, we clearly have $\el - \epsilon \ge 0$. 
In either case, the inequality $\el - \epsilon \ge 0$ holds. 
Then, 
for a vertex $v$ of $V(D^{*}) \setminus \{t^{*}, t_{c\hspace*{-0.25mm}d}, z\}$, 
the following hold:
\begin{itemize}
\item 
If $ |N_{G}(v) \cap T|  = 0$, 
then 
$\deg_{D^{*}}(v) = \deg_{G}(v) \ge \deg_{G}(v) - \el + \epsilon$.
\item 
If $1 \le |N_{G}(v) \cap T| \le \el$, 
then 
$\deg_{D^{*}}(v) \ge \deg_{G}(v) - \el+1 \ge \deg_{G}(v)- \el +\epsilon$. 
\item
If $|N_{G}(v) \cap T| = \el+1$\footnote{In this case, if $\epsilon = 1$ then  
$v$ is adjacent to all the vertices of $T$.},
then 
$\deg_{D^{*}}(v) \ge \deg_{G}(v) - (\el+1) + (1+ \epsilon) = \deg_{G}(v) - \el +\epsilon $.
\end{itemize} 
%
%
%
%
%
%
\noindent
Thus 
the definition of $D^{*}$ and Claim~\ref{claim:at most l+1} yield that
\[
\delta(D^*, t^{*}, t_{c\hspace*{-0.25mm}d}; z) 
\ge 
\delta(G, x, y; z)-\el + \epsilon \ge (k-\el+\epsilon)+1. 
\]
By the induction hypothesis, 
$D^*$ contains $k - \el + \epsilon$ admissible $(t^{*}, t_{c\hspace*{-0.25mm}d})$-paths. 
This implies that 
$G[T \cup D]$ contains $k-\el+\epsilon$ admissible $(T \setminus \{t_{c\hspace*{-0.25mm}d}\}, t_{c\hspace*{-0.25mm}d})$-paths $\oraP{1}, \dots, \oraP{k-\el+\epsilon}$. 
Let 
$t_{i}$ be the unique vertex of $V(P_{i}) \cap (T \setminus \{t_{c\hspace*{-0.25mm}d}\})$ for $1 \le i \le k-\el+\epsilon$. 
Then  
$t_{i}\oraP{i}t_{c\hspace*{-0.25mm}d}y$
($1 \le i \le k-\el+\epsilon$) are $k-\el+\epsilon$ admissible $(T \setminus \{t_{c\hspace*{-0.25mm}d}\}, y)$-paths in $G[T \cup D \cup C]$. 
On the other hand, 
it follows from Fact~\ref{fact:paths in H}~\ref{xT} that 
for each $1 \le i \le k-\el+\epsilon$, 
$H - t_{c\hspace*{-0.25mm}d}$ contains $\el-\epsilon+1$ $(x, t_{i})$-paths
of 
lengths $1, 2, \ldots, \el-\epsilon+1$ (if $H$ is of type~1)
or $1,3, \ldots, 2(\el-\epsilon)+1$ (if $H$ is of type~3).
Hence by Fact~\ref{fact:Lemma3.2}, 
we obtain 
$k \ \big( = (k-\el+\epsilon) + (\el-\epsilon+1) - 1 \big)$ admissible $(x, y)$-paths in $G$,
a contradiction. 
Thus $N_{G}(D) \not\subseteq T$. 
Combining this with Claim \ref{claim:neighbor of x and y}, 
we have
$N_{G}(D) \cap S 
= N_{G}(D) \setminus (T \cup x) 
= N_{G}(D) \setminus T 
\not=\emptyset$.
\qed

\medskip
\noindent\textbf{Case~1.1.}
$H$ is of type~1.

By Claim~\ref{claim:neighbor of x and y}, 
$N_{G}(x) = N_{G}(y) = T$. 
By Claim~\ref{claim:notsubsetT}, 
we also have 
$V(G) = T \cup \{x, y, z\}$.
Since 
$|T| = \el + 1 \le k - 1$ by Fact~\ref{fact:l >= k or l >= k-1}~(2), 
and since $T \sm \{z\} \neq \emptyset$, 
the degree condition yields that 
$(\el + 1 = ) \ |T| = k - 1$, 
$z \notin T \cup \{x, y\}$, 
and 
$N_{G}(v) = (T \setminus \{v\}) \cup \{x, y, z\}$ for all $v \in T$. 
This implies that
$G$ contains $(x, y)$-paths of lengths $2,3,\ldots,k + 1$. 
Thus $G$ contains $k$ admissible $(x, y)$-paths, a contradiction.

\medskip
\noindent\textbf{Case~1.2.}
$H$ is of type~2.

By Claim~\ref{claim:neighbor of x and y}, 
we have $N_{G}(x) = N_{G}(y) = S$, 
say $N_{G}(x) = N_{G}(y) = \{s_{1}, s_{2}\}$. 
Let $G' = G - \{x, y\}$. 
Since $G$ and $H-x$ are $2$-connected, respectively, 
and 
$|V(H-x)| \ge 4$, 
it follows that 
$(G', s_{1}, s_{2}; z)$ is a $2$-connected rooted graph 
such that 
$\delta(G', s_{1}, s_{2}; z) \ge \delta(G, x, y; z)$.
Therefore, 
by the induction hypothesis, 
we obtain $k$ admissible $(s_{1}, s_{2})$-paths $\oraP{1}, \dots, \oraP{k}$ in $G'$.
Then 
$xs_{1}\oraP{i}s_{2}y$ ($1 \le i\le k$)
are $k$ admissible $(x, y)$-paths in $G$,
a contradiction.

\medskip
\noindent\textbf{Case~1.3.}
$H$ is of type~3.

\begin{claim}\label{claim:D}
There exists a component $D$ of $G - V(H)$
such that $D\not=C$, $V(D) \setminus \{z\} \neq \emptyset$ 
and $N_{G}(D) \cap T \not= \emptyset$. 
\end{claim}

\proof 
If there exists $t \in T$ such that 
$N_{G}(t) \setminus (V(H) \cup \{y, z\}) \neq \emptyset$, 
then the assertion clearly holds. 
Thus, 
we may assume that 
$N_{G}(t) \setminus (V(H) \cup \{y, z\}) = \emptyset$ 
for all $t \in T$. 
Since 
$N_{G}(y) \cap T \neq \emptyset$ by Claim~\ref{claim:neighbor of x and y},
Fact~\ref{fact:l >= k or l >= k-1}~(2) 
yields $|S| = \el \le k-2$. 
Then 
for a vertex $t \in T \setminus \{z\} \ (\neq \emptyset)$, 
we have 
\[
0 
= |N_{G}(t) \setminus (V(H) \cup \{y, z\})| 
\ge 
(k + 1) - (|S| + |\{x, y, z\}|) 
\ge (k + 1) - (k + 1) 
= 0.
\]
Thus the equality holds, 
which implies that 
$|S| = \el = k - 2$, $z \notin V(H) \cup \{y\}$ 
and 
$tz \in E(G)$ for all $t \in T$. 
By Claim~\ref{claim:neighbor of x and y} 
and 
since there exists no core of type~2 with respect to $(x, y)$ by \ref{renfomin}, 
we also have $N_{G}(x)=N_{G}(y)=T = N_{G}(z) \cap V(H)$. 

If $N_{G}(z) \sm V(H) \not= \emptyset$,
then since $T \subseteq N_{G}(z)$, 
the claim follows, 
and so we may assume that 
$N_{G}(z) \sm V(H) = \emptyset$,
that is, $N_{G}(z) = T$. 
If $|T| \ge \el+2$, 
then
$x \vee T \vee (S \cup z)$
is an $(\el + 1)$-core of type~3, 
contradicting to \ref{TSlarge}-(i).
Thus we have $|T| = \el + 1 = k - 1$, 
which also implies that $|S| = \el \ge 1$. 
Since 
$N_{G}(x)=N_{G}(y)=N_{G}(z) = T$, 
a vertex $s \in S$ satisfies 
\[
|N_{G}(s) \setminus (V(H) \cup \{y, z\})|
= 
|N_{G}(s) \setminus T|
\ge 
(k + 1) - |T| = (k + 1) - (k-1)
> 0. 
\]
Hence there exists a component $D$ of $G - V(H)$
such that $D\not=C$, $V(D) \setminus \{z\} \neq \emptyset$,
and $N_{G}(D) \subseteq S$. 
Let $s_{0} \in N_{G}(D)$ and
$D^{*}$ be the graph obtained from $G[D \cup N_{G}(D)]$
by contracting $N_{G}(D) \sm \{s_{0}\}$
into a new vertex $s^{*}$. 
Since $N_{G}(D) \subseteq S$ and 
$G$ is $2$-connected, 
it follows that $(D^{*}, s^{*}, s_{0} ; z)$ is a $2$-connected rooted graph. 
Since $V(D) \setminus \{z\} \neq \emptyset$ 
and $|S| = \el = k-2$, 
we also have 
$\delta(D^{*}, s^{*}, s_{0}  ; z) 
\ge 
\delta(G, x, y; z) - (k-2) + 1 \ge 3 +1$. 
Therefore, by the induction hypothesis, 
$D^{*}$ contains three admissible $(s^{*}, s_{0})$-paths. 
This implies that 
$G[S \cup D]$ contains three admissible $(S \setminus \{s_{0}\}, s_{0})$-paths $\oraP{1}, \oraP{2}$ and $\oraP{3}$. 
Let 
$s_{i}$ be the unique vertex of $V(P_{i}) \cap (S \setminus \{s_{0}\})$ for $1 \le i \le 3$, 
and let $t_{0} \in N_{G}(y) \cap T$. 
Then  
$s_{i}\oraP{i}s_{0}t_{0}y$ ($1 \le i \le 3$) are
 three admissible $(S \setminus \{s_{0}\}, y)$-paths in $G[t_{0} \cup S \cup D \cup C]$. 
On the other hand, 
since $|T \sm \{t_{0}\}|=| (S \sm \{s_{0}\}) \cup \{z\}| = \el = k - 2$ 
and 
$N_{G}(z) = T$, 
it follows that 
for each $1 \le i \le 3$, 
$G[\big( V(H) \setminus \{t_{0},s_{0}\} \big)\cup \{z\}]$ contains $k - 2$ $(x, s_{i})$-paths
of lengths $2, 4, \ldots, 2(k-2)$.
Hence by Fact~\ref{fact:Lemma3.2}, 
we obtain 
$k \ \big( = 3 + (k-2) - 1 \big)$ admissible $(x, y)$-paths in $G$,
a contradiction. 
%
%
%
%
%
%
\qed

Let 
$D$ be a component of $G - V(H)$ as in Claim~\ref{claim:D}. 
Since 
\ref{transform} does not hold, 
it follows from Claims~\ref{claim:neighbor of x and y} and \ref{claim:D} that 
$|T|=2$, say $T = \{t_{1}, t_{2}\}$. 
Since 
$N_{G}(D) \cap S\not=\emptyset$ by Claim~\ref{claim:notsubsetT} 
and since 
$|T| \ge |S|+1$, 
we also have $|S| = 1$, say $S = \{s\}$. 
Let $G' = G - \{x, y\}$. 
Since $G$ is $2$-connected 
and $N_{G}(x) = N_{G}(y) = \{t_{1}, t_{2}\}$ by Claim~\ref{claim:neighbor of x and y}, 
it is easy to check that $(G', t_{1}, t_{2}; z)$ is a $2$-connected rooted graph. 
Since 
$\emptyset \neq V(D) \setminus \{z\} \subseteq V(G')$, 
we also have 
$\delta(G', t_{1}, t_{2}; z) \ge \delta(G, x, y; z)$. 
Therefore, 
by the induction hypothesis, 
we obtain $k$ admissible $(t_{1}, t_{2})$-paths $\oraP{1}, \dots, \oraP{k}$ in $G'$.
Then 
$xt_{1}\oraP{i}t_{2}y$ ($1 \le i\le k$)
are $k$ admissible $(x, y)$-paths in $G$,
a contradiction.

This completes the proof of Case~1.

\bigskip
\noindent\textbf{Case~2.}
$V(\CF) \not= \{y\}$.

\begin{claim}\label{claim:S=1NCT}
Assume that $H$ is of type~3. 
If $|S| = 1$, then $N_{G}(\CF) \cap \TF \neq \emptyset$.
\end{claim}
\proof 
Suppose that $|S| = 1$, say $S = \{s\}$, 
and $N_{G}(C) \cap T = \emptyset$. 
Let $G' = G - V(C)$. 
Since 
$G$ and $H$ are $2$-connected, 
$y \notin V(G')$ 
and 
$|V(H)| \ge |\{x\}| + |T| + |S| \ge 1 + 2 + 1 \ge 4$, 
it follows that 
$(G', x, s; z)$ is a $2$-connected rooted graph 
such that 
$\delta(G', x, s;z) \ge \delta(G, x, y; z)$. 
By the induction hypothesis, 
$G'$ contains $k$ admissible $(x, s)$-paths $\oraP{1}, \ldots, \oraP{k}$. 
Since 
$G$ is $2$-connected 
and $N_{G}(C) \cap T = \emptyset$, 
we have 
$s \in N_{G}(C)$, 
and so there exists an $(s,y)$-path $\ora{Q}$ in $G[C \cup s]$. 
Then
$x\oraP{i}s\ora{Q}y$ ($1 \le i \le k$)
are $k$ admissible $(x, y)$-paths in $G$, a contradiction. 
\qed

In this case, 
we will apply the induction hypothesis 
for new graphs obtained 
from $H$ and 
blocks with at most two cut-vertices of $\CF$. 
However, 
the $z$-end-block of $\CF$ will not help us to find admissible paths in the argument.  
So, 
in the following two claims, 
we study the structure for the case where $\CF$ contains the $z$-end-block. 
In particular, we show that 
$\CF$ is not a $(y,z)$-path of order exactly $3$
at this stage. 
(See Subsection~\ref{subsec:terminology and notation} 
for the definitions of the $z$-end-block $B_{z}$ and the vertices $b_{z}, b_{z}'$.)

\begin{claim}\label{claim:notK2}
Assume that 
there exists the $z$-end-block $B_{z}$ with cut-vertex $b_{z}$ in $\CF$ such that $y \notin \{z, b_{z}\}$. 
Assume further that 
$\deg_{\CF}(b_{z}) = 2$. 
Then the following hold. 
\begin{center}
\begin{tabular}{p{65mm}p{65mm}}
(1) $\ell = k - 2$.
&
(2) $|N_{G}(b_{z}) \cap V(\HF)| = \el+1$. \\[3mm]
(3) $\big( N_{G}(z) \cup N_{G}(b_{z}') \big) \cap \TF = \emptyset$.
&
(4) If $b_{z}' \not= y$,
then 
$\deg_{\CF}(b_{z}') \ge 3$.
\end{tabular}
\end{center}
\end{claim}
\proof 
By our assumption, 
$\deg_{G}(b_{z}) \ge k + 1$ 
and 
there exists a $(b_{z},y)$-path $\ora{R}$ in $\CF - z$. 
If $H$ is of type 3$^{\fl}$, 
then since 
$y\not=z$, $N_{\CF}(y) = \{t_{0}\}$ and 
by Claim~\ref{claim:transformed core}, 
note that $b_{z} \neq t_{0}$. 

\smallskip
(1),(2) To show (1) and (2), we first prove that 
\begin{align}
\label{el <= k-2}
\el \le k - 2. 
\end{align}
Since 
$\el \le k - 1$ by Fact~\ref{fact:l >= k or l >= k-1}~(1), 
it suffices to show that $\el \neq k - 1$. 
Suppose to the contrary that $\el = k - 1$. 
Then 
it follows from Fact~\ref{fact:l >= k or l >= k-1}~(2) that $N_{G}(C) \cap T = \emptyset$. 
Combining this with 
Claim~\ref{claim:2-connected and xy notin E}, 
we have $N_{G}(z) \cap S \not=\emptyset$, 
say $s_{z} \in N_{G}(z) \cap S$. 
This in particular implies that $H$ is of type~2 or type~3. 

Suppose that $N_{G}(b_{z}) \cap S \not=\emptyset$, 
say $s_{b} \in N_{G}(b_{z}) \cap S$. 
Then 
$s_{b}b_{z}\ora{R}y$ and $s_{z}zb_{z}\ora{R}y$ are two admissible $(\{s_{b}, s_{z}\}, y)$-paths in $G[S \cup C]$.
On the other hand, 
it follows from 
Fact~\ref{fact:paths in H}~(1) that 
for each $s \in \{s_{b}, s_{z}\}$, 
$H$ contains $k-1 \ (= \el)$ admissible $(x, s)$-paths. 
Hence by Fact~\ref{fact:Lemma3.2}, 
we obtain $k \ \big( = 2 + (k - 1) -1\big)$ admissible $(x, y)$-paths in $G$,
a contradiction. 
Thus $N_{G}(b_{z}) \cap S =\emptyset$, 
that is, $N_{G}(b_{z}) \cap (S \cup T) =\emptyset$. 
Then 
$1 \le k - 1 \le \deg_{G}(b_{z}) - \deg_{C}(b_{z}) = |N_{G}(b_{z}) \cap V(H)| \le |\{x\}| = 1$. 
Thus the equality holds, 
which implies that 
$\el = k - 1 = 1$ and $N_{G}(b_{z}) \cap V(H) = \{x\}$. 
If $H$ is of type~2, 
then 
$xb_{z} \ora{R}y$ 
and 
$x s_{z} z b_{z} \ora{R}y$ 
are $k \ (= 2)$ admissible $(x, y)$-paths in $G$, a contradiction; 
if $H$ is of type~3, 
then since 
$|S| = \el = 1$ 
and $N_{G}(C) \cap T = \emptyset$, 
this contradicts Claim~\ref{claim:S=1NCT}.
Thus (\ref{el <= k-2}) is proved.  

Now, by Claim~\ref{claim:at most l+1} 
and (\ref{el <= k-2}), 
we have 
$$k - 1 \le \deg_{G}(b_{z}) - \deg_{\CF}(b_{z}) = |N_{G}(b_{z}) \cap V(\HF)| \le \ell + 1 \le k-1.$$
Thus the equality holds, which implies that 
$\el = k-2$ and $|N_{G}(b_{z}) \cap V(\HF)| = \el + 1$.

\smallskip
(3) 
Note that by Claims~\ref{claim:at most l+1}
and \ref{claim:notK2}~(2), 
$N_{G}(b_{z}) \cap \TF \neq \emptyset$, 
say $t_{b} \in N_{G}(b_{z}) \cap \TF$. 
To show (3), 
suppose that $N_{G}(v) \cap \TF \neq \emptyset$ for some $v \in \{z, b_{z}'\}$, 
and let $t_{v} \in N_{G}(v) \cap \TF$. 

Since $N_{\CF}(b_{z}) = \{z, b_{z}'\}$, 
it follows that $G[\{z, b_{z}, b_{z}', t_{v}\}]$ contains a $(t_{v}, b_{z}')$-path $\ora{P}$ of length $1$ or $3$. 
Hence $P$ and $t_{b} b_{z} b_{z}'$ are $(\{t_{v}, t_{b}\}, b_{z}')$-paths of lengths $1, 2$ or $3, 2$. 
By adding $b_{z}' \ora{R} y$ to each of the two paths, 
we obtain two semi-admissible $(\{t_{v}, t_{b}\}, y)$-paths in $G[\CF \cup \{t_{v}, t_{b}\}]$. 
On the other hand, 
it follows from 
Fact~\ref{fact:paths in H}~(2) and Claim~\ref{claim:notK2}~(1) that 
for each $t \in \{t_{v}, t_{b}\}$, 
$\HF$ contains $k-1 \ (= \el + 1)$ semi-admissible $(x, t)$-paths. 
Hence by Fact~\ref{fact:Lemma3.2}, 
$G$ contains $k \ \big( = 2 + (k-1) - 1 \big)$ admissible $(x, y)$-paths, a contradiction.

\smallskip
(4) Assume that $b_{z}' \not= y$
and $\deg_{\CF}(b_{z}') \le 2$. 
We first claim that $b_{z}' \neq t_{0}$ if $H$ is of type~$3^{\fl}$. 
Suppose to the contrary that 
$H$ is of type~$3^{\fl}$ and $b_{z}'=t_{0}$.
(See Figure \ref{transformFIG}.)
If $H$ is of type~$3^{\fl}$ in \ref{M1}, 
then 
$\deg_{C}(b_{z}') = \deg_{C}(t_{0}) 
\ge 
|N_{G}(t_{0}) \cap V(D_{0})|
+ |\{y, s_{0}\}| \ge 1 + 2 = 3$, 
a contradiction.
Thus $H$ is of type~$3^{\fl}$ in \ref{M2}.
Recall that $(b_{z}b_{z}' = ) \ b_{z}t_{0} \in E(G)$ 
and 
$x \cup S \subseteq N_{G}(t_{0})$. 
Since 
there exist no cores of type 1 or type~2 with respect to $(x, y)$ by \ref{renfomin}, 
this together with Claim~\ref{claim:notK2}~(2) 
implies that 
$|N_{G}(b_{z}) \cap T| = |N_{G}(b_{z}) \cap V(H)| = \el + 1$. 
Hence 
$H' 
:= 
x \vee \big( (N_{G}(b_{z}) \cap T) \cup t_{0} \big) \vee (S \cup b_{z})$
is an $(\el + 1)$-core of type~3, 
which contradicts \ref{TSlarge}-(i). 
Thus $b_{z}'\neq t_{0}$ if $H$ is of type~$3^{\fl}$.

By Claim~\ref{claim:notK2}~(3), $N_{G}(b_{z}') \cap \TF = \emptyset$, 
and so Claims~\ref{claim:at most l+1}  and \ref{claim:notK2}~(1)
yield that $|N_{G}(b_{z}') \cap V(\HF)| \le \el=k-2$. 
Then we obtain
\[
\deg_{G}(b_{z}') \le \deg_{C}(b_{z}') + |N_{G}(b_{z}') \cap V(\HF)|
\le 2+\el= k ,
\]
a contradiction. 

\indent
This completes the proof of Claim~\ref{claim:notK2}. 
\qed

\begin{claim}\label{claim:notYZ}
$\CF$ is not a $(y,z)$-path of order exactly $3$. 
\end{claim}
\proof 
Suppose that $C$ is a $(y,z)$-path of order exactly $3$. 
By 
Claims~\ref{claim:2-connected and xy notin E} 
and \ref{claim:notK2}~(3), 
we have $N_{G}(z) \cap S \neq \emptyset$, say $s_{z} \in N_{G}(z) \cap S$. 
This in particular implies that $H$ is not of type~1.

Suppose that $H$ is of type~2. 
Since $N_{G}(b_{z}) \cap T \neq \emptyset$ by Claims~\ref{claim:at most l+1} and \ref{claim:notK2}~(2), 
it follows from Fact~\ref{fact:paths in H}~(2) and Claim~\ref{claim:notK2}~(1) that 
$G[H \cup b_{z}]$ contains $k-1 \ (= \el + 1)$ admissible $(x, b_{z})$-paths 
$\oraP{1}, \dots, \oraP{k-1}$ of lengths $3,4,\ldots,\el+3$. 
On the other hand, 
by Fact \ref{fact:paths in H}~(1),
$H$ contains an $(x, s_{z})$-path $\ora{Q}$ of length $\el+2$, 
and so 
$\oraP{k} := x\ora{Q}s_{z}zb_{z}$ 
is an $(x, b_{z})$-path of length $\el+4$. 
Then
$x\oraP{i}b_{z}y$ ($1 \le i \le k$)
are $k$ admissible $(x, y)$-paths in $G$, 
a contradiction. 
Thus $H$ is not of type~2, that is, $H$ is of type~3.

Since $N_{G}(y) \cap T \ (= N_{G}(b_{z}') \cap T)=\emptyset$ by Claim~\ref{claim:notK2}~(3), 
and since $xy \notin E(G)$ by Claim \ref{claim:2-connected and xy notin E}~(2), 
it follows that 
$N_{G}(y) \subseteq S \cup b_{z}$. 
By \ref{refno} and \ref{renfomin}, 
there exist no cores of type 1 or type~2 with respect to $(y, x)$, 
and so 
any core with respect to $(y, x)$ is of type $3$ (by Remark~\ref{remark:lcore}, 
Claim~\ref{claim:2-connected and xy notin E}). 
Then we can use the inequality in \ref{degxy}. 
Note that $\el \ge 1$, since $S \neq \emptyset$. 
Hence 
\[
\el + 1 = \max\{\el + 1, 2\} \le |T| \le \deg_{G}(x) \le \deg_{G}(y) \le |S \cup b_{z}| = \el + 1. 
\]
Thus the equality holds, 
which implies that $\deg_{G}(x) = \deg_{G}(y)$ and $N_{G}(x) = T$. 
Then it follows from the first equality and \ref{dist} that 
$\dist_{G}(x, z) \le \dist_{G}(y, z) =2$ holds. 
On the other hand, 
since $xz \notin E(G)$ by Claim~\ref{claim:2-connected and xy notin E}~(2),
and since 
$N_{G}(z) \cap N_{G}(x) = N_{G}(z) \cap T = \emptyset$ 
by Claim~\ref{claim:notK2}~(3), 
it follows that 
$\dist_{G}(x, z) \ge 3$. 
This is a contradiction. 
\qed

Let $\CT$ be the set of cut-vertices of $C$. 
A block $B$ of $\CF$ is said to be \textit{feasible} 
if $B$ satisfies the following condition~\ref{FEB}. 
\begin{enumerate}[label=(\texttt{F})]
\item
\label{FEB}
$|V(B) \cap (\CT \cup \{y, z\})| \le 2$ 
and 
$V(B) \setminus (\CT \cup \{y, z\}) \neq \emptyset$. 
\end{enumerate}
Note that by the assumption of Case~2 and Claim~\ref{claim:2-connected and xy notin E}~(2), 
if $\CF$ itself is a 
block, then $\CF$ is feasible.

\begin{claim}\label{claim:feasible}
There exists 
a feasible block
of $\CF$.
\end{claim}

\proof
Suppose that there exists no feasible block of $\CF$. 
Then the condition~\ref{FEB} yields the following:
$\CF$ is not a block; 
its block-tree is a path; 
one of the two end-blocks of $\CF$ is the $y$-end-block 
and the other is the $z$-end-block. 
By the definition of $b_{z}$ and $b_{z}'$, 
if $\deg_{\CF}(b_{z}) \ge 3$, 
then 
$b_{z} = b_{z}'$ 
and so 
$\deg_{\CF}(b_{z}') \ge 3$; 
if 
$\deg_{\CF}(b_{z}) = 2$, 
then 
it follows from Claims~\ref{claim:notK2}~(4) 
and  \ref{claim:notYZ}
that $\deg_{\CF}(b_{z}') \ge 3$. 
In either case, 
$\deg_{\CF}(b_{z}') \ge 3$ holds. 
Hence 
there exists a block $B$ of $C$ 
which is not an end-block 
and satisfies \ref{FEB}.
\qed

\begin{figure}[H]
\begin{center}
{\includegraphics[pagebox=artbox]{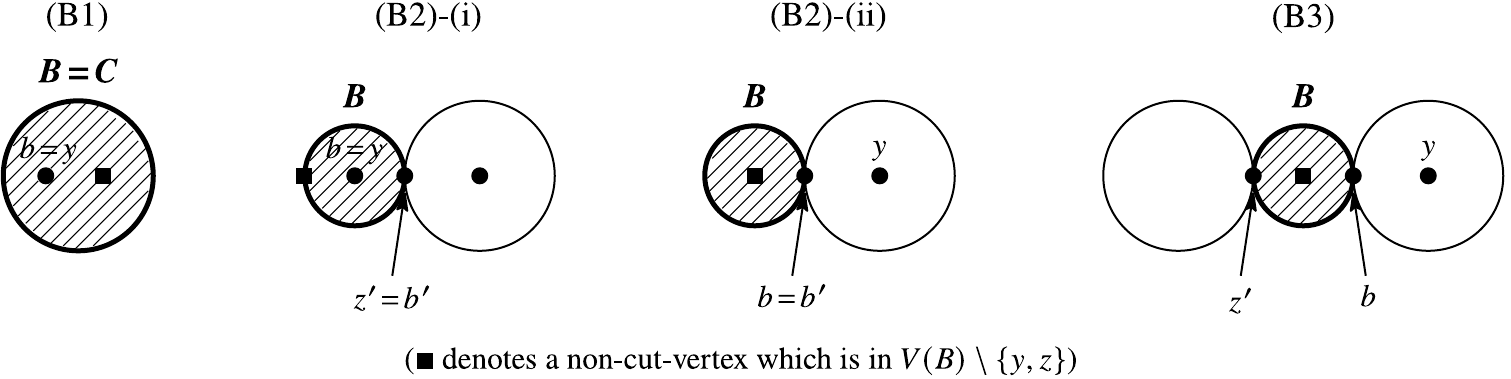}}
\end{center}
\caption{The definitions of $B, b$ and $z'$}
\label{feasibleFIG}
\end{figure}

In the rest of the proof, 
$B$, $b$ and $z'$ denote any one of the following~\ref{Cis2-conn}, \ref{existsFEB}-\ref{yinB}, 
\ref{existsFEB}-\ref{ynotinB} and \ref{existsFB} 
(note that by Claim~\ref{claim:feasible} and \ref{FEB}, such a tuple $(B, b, z')$ exists, 
see also Figure~\ref{feasibleFIG}):
\begin{enumerate}[label=(\texttt{B}\arabic*)]
\item
\label{Cis2-conn}
$B$ is a feasible block of $\CF$ such that $|V(B) \cap \CT| = 0$ 
(i.e., $\CF$ itself is a block and $B = C$)
and, $b := y$ and $z' := z$. 
\item
\label{existsFEB}
$B$ is a feasible block of $\CF$ such that $|V(B) \cap \CT| = 1$, say $V(B) \cap \CT = \{b'\}$,  
and 
\begin{enumerate}[label=(\roman*)]
\item
\label{yinB}
if $y \in V(B) \setminus \{b'\}$,
then
$b := y$ 
and 
$z' := b'$; 
\item
\label{ynotinB}
if $y \not\in V(B) \setminus \{b'\}$,
then
$b := b'$ 
and 
$z' :=z$. 
\end{enumerate}
\item
\label{existsFB}
$B$ is a feasible block of $\CF$ such that $|V(B) \cap \CT| = 2$, 
and 
$b$ is the unique vertex of $V(B) \cap \CT$ 
such that 
$C - ( V(B) \setminus \CT)$ 
contains a $(b, y)$-path 
(possibly $b = y$) 
and 
$\{z'\} := (V(B) \cap \CT) \setminus \{b\}$. 
\end{enumerate}
Note that 
$\emptyset \neq V(B) \setminus \{b, z'\} \subseteq V(G) \setminus \{x, y, z\}$ 
and 
$N_{G}(v) \subseteq V(B) \cup V(H)$ for $v \in V(B) \setminus \{b, z'\}$. 
Note also that 
if $\HF$ is of type~$3^{\fl}$, 
then since $t_{0}$ is a cut-vertex of $\CF$, 
$t_{0} \notin V(B) \setminus \{b, z'\}$.
Let $\ora{R}$ be a $(b,y)$-path in $C$ such that $V(R) \cap V(B - b) = \emptyset$.

\begin{claim}\label{claim:subsetSbx}
(1) $N_{G}(B-b) \cap V(\HF) \subseteq x \cup S$
and
(2) $|N_{G}(B-b) \cap (x \cup S)| \ge 2$.
\end{claim}

\proof 
(1) Suppose that $N_{G}(B-b) \cap \TF \not=\emptyset$.
Let $B^{*}$ be the graph obtained from 
$G[B \cup \big( N_{G}(B-b) \cap \TF \big)]$
by contracting $N_{G}(B-b) \cap \TF$ into a new vertex $t^{*}$.
Since $\emptyset \neq V(B) \setminus \{b, z'\} \subseteq V(G) \setminus \{x, y, z\}$, 
$(B^{*}, t^{*}, b ; z')$ is a $2$-connected rooted graph 
such that $\emptyset \neq V(B^{*}) \setminus \{t^{*}, b, z' \} \subseteq V(G) \setminus \{x, y, z\}$.
Then, it follows from Claim~\ref{claim:at most l+1} that 
for a vertex $v$ of $V(B^{*}) \setminus \{ t^{*}, b , z'\}$, 
the following hold:
\begin{itemize}
\item 
If $ |N_{G}(v) \cap \TF|  = 0$, 
then 
$\deg_{B^{*}}(v) \ge \deg_{G}(v)-\el$.
\item 
If $|N_{G}(v) \cap \TF| \ge 1$, 
then 
$\deg_{B^{*}}(v) \ge \deg_{G}(v) - (\el + 1) +1 = \deg_{G}(v) - \el$.
\end{itemize}
Thus 
the definition of $B^{*}$ yields that
\[
\delta(B^{*}, t^{*}, b ; z') 
\ge 
\delta(G, x, y; z) - \el  \ge (k - \el)+1. 
\]
By the induction hypothesis, 
$B^{*}$ contains $k - \el$ admissible $(t^{*}, b)$-paths.
Thus 
$G$ contains $k - \el$ admissible $(\TF, b)$-paths
$\oraP{1}, \dots, \oraP{k - \el}$. 
Let $t_{i} \in V(P_{i}) \cap \TF$ for $1 \le i \le k - \el$. 
Then 
$t_{i} \oraP{i} b \ora{R} y$ is a $(t_{i}, y)$-path in $G[\CF \cup t_{i}]$ for $1 \le i \le k - \el$. 
On the other hand, 
by Fact~\ref{fact:paths in H}~\ref{xT},
$\HF$ contains $\el+1$ semi-admissible $(x, t_{i})$-paths for $1 \le i \le k - \el$.
Hence by Fact~\ref{fact:Lemma3.2}, 
$G$ contains 
$k \ \big( = (k - \el) + (\el + 1) - 1 \big)$ admissible $(x, y)$-paths, a contradiction. 
Thus (1) holds.

(2) Suppose that 
either (i) $|N_{G}(B-b) \cap (x \cup S)| = 1$ 
or 
(ii) $N_{G}(B-b) \cap (x \cup S) = \emptyset$ holds. 
Note that if $B$ satisfies \ref{Cis2-conn} or \ref{existsFEB}-\ref{ynotinB}, 
then the $2$-connectivity of $G$ implies that (i) holds; 
that is to say, 
if (ii) holds, then 
$B$ satisfies \ref{existsFEB}-\ref{yinB} or \ref{existsFB}. 
If (i) holds, 
say $N_{G}(B-b) \cap (x \cup S)=\{v\}$, 
then let $B' = G[B \cup v]$; 
if (ii) holds, 
then let $v = z'$ and $B' = B$. 
Since $\emptyset \neq V(B) \setminus \{b, z'\} \subseteq V(G) \setminus \{x, y, z\}$, 
$(B', v, b; z')$ is a $2$-connected rooted graph such that 
$\delta(B', v, b; z') \ge \delta(G, x, y; z)$. 
By the induction hypothesis, 
$B'$ 
contains $k$ admissible $(v, b)$-paths 
$\oraP{1}, \dots, \oraP{k}$. 
If (i) holds, then 
let $\ora{Q}$ be an $(x, v)$-path in $\HF$; 
if (ii) holds, then 
by the $2$-connectivity of $G$, 
there exists an $(x, v)$-path $\ora{Q}$ 
in 
$G[\HF \cup \big( V(\CF) \setminus ( V(B' - v) \cup V(R) ) \big) ]$. 
Then
$x\ora{Q}v\oraP{i}b\ora{R}y$ ($1 \le i \le k$)
are $k$ admissible $(x, y)$-paths in $G$, a contradiction. 
\qed

\medskip
\noindent\textbf{Case~2.1.}
$H$ is of type~1.

By Claim~\ref{claim:subsetSbx}~(2), 
we have
$N_{G}(B-b) \cap S \not= \emptyset$,
which 
contradicts $S=\emptyset$.

\medskip
\noindent\textbf{Case~2.2.}
$H$ is of type~2.

By Claim~\ref{claim:subsetSbx}~(2), 
$N_{G}(B-b) \cap S \not= \emptyset$. 
Let $B^{*}$ be the graph obtained from 
$G[B \cup (N_{G}(B-b) \cap S)]$
by contracting $N_{G}(B-b) \cap S$ into a new vertex $s^{*}$. 
Then 
$(B^{*}, s^{*}, b; z')$ is a $2$-connected rooted graph 
such that $V(B^{*}) \setminus \{s^{*}, b, z'\} = V(B) \setminus \{b, z'\} \neq \emptyset$. 
Since 
there exists no core of type~1 with respect to $(x, y)$ by \ref{renfomin}, 
it follows that 
$x \not\in N_{G}(v)$ or $N_{G}(v) \cap S =\emptyset$ for $v \in V(B) \setminus \{b\}$. 
This together with the definition of $B^{*}$ and Claim~\ref{claim:subsetSbx}~(1) 
implies that 
$\delta(B^{*}, s^{*}, b; z') \ge \delta(G, x, y; z) - 1 \ge (k - 1) + 1$. 
Hence, 
by the induction hypothesis, 
$B^{*}$ contains $k-1$ admissible $(s^{*}, b)$-paths, 
and so 
$G[S \cup B]$ contains 
$k -1$ admissible $(S, b)$-paths $\oraP{1}, \dots, \oraP{k-1}$. 
Let $s_{i} \in V(P_{i}) \cap S$ for $1 \le i \le k-1$. 
Then $s_{i} \oraP{i} b \ora{R} y$ is an $(s_{i}, y)$-path in $G[C \cup s_{i}]$ for $1 \le i \le k-1$. 
On the other hand, 
by Fact~\ref{fact:paths in H}~\ref{xS}, 
$H$ contains two admissible $(x, s_{i})$-paths 
for $1 \le i \le k - 1$. 
By Fact~\ref{fact:Lemma3.2}, 
$G$ contains $k \ \big( = (k-1) + 2 - 1 \big)$ admissible $(x, y)$-paths, a contradiction.

\medskip
\noindent\textbf{Case~2.3.}
$H$ is of type~3.

Note that, 
by Claim~\ref{claim:subsetSbx}~(2), 
$\el = |S| \ge 1$. 
Let 
\[
\NCT = V(\CF) \setminus (\CT \cup \{y, z\}).
\]
We divide $\HF$ into three cases as follows:
\begin{itemize} 
\item
$\HF$ \textit{is of type~I} 
if 
$|N_{G}(v_{0}) \cap \TF| =\el+1$ 
for some 
$v_{0} \in \NCT$. 
\item
$\HF$ \textit{is of type~II} 
if 
$\el=1$ and 
$|N_{G}(v_{0}) \cap S| =|N_{G}(v_{0}) \cap \TF|=1$ 
for some 
$v_{0} \in \NCT$. 
\item
$\HF$ \textit{is of type~III}
if $H$ is of neither type~I nor type~II. 
\end{itemize} 
%
%
%
%
%
%
If $\HF$ is of type~I or type~II, 
then let $v_{0}$ be a vertex as described above, and let $S^{\sh} =S \cup v_{0}$;
if $\HF$ is of type~III, then let $S^{\sh} = S$. 
We then let 
\[\text{
$H^{\sh} = G[x \cup \TF \cup S^{\sh}]$ \, and \, $\lp =|S^{\sh}|$. 
}\]
Then the following (i) and (ii) hold: 
(i) $\lp \ge \el \ge 1$, 
and 
if $\HF$ is of type~I or type~II, then $\lp \ge 2$; 
(ii) if $\HF$ is of type~I or type~II, 
then by the definitions of $\NCT$ and the types, 
and by Claim~\ref{claim:subsetSbx}~(1), 
$v_{0} \notin V(B)$ 
and 
$v_{0}$ does not separate $B$ and $y$ in $C$. 
In particular, by (ii), 
$B$ is still a block of a component of $G - V(H^{\sh})$, 
there exists a $(b, y)$-path 
internally disjoint from $B$ 
in $G - V(H^{\sh})$, 
and 
$N_{G}(v) \subseteq V(B) \cup (x \cup S)$ for $v \in V(B) \setminus \{b, z'\}$ 
(by Claim~\ref{claim:subsetSbx}~(1)).

\begin{claim}\label{claim:l=1}
$\lp =1$.
\end{claim}

\proof 
Suppose that $\lp \ge 2$. 
Let $B^{*}$ be the graph obtained from 
$G[B \cup \big(N_{G}(B-b) \cap S \big)]$
by contracting $N_{G}(B-b) \cap S$ into a new vertex $s^{*}$.
By Claim~\ref{claim:subsetSbx}~(2),
$(B^{*}, s^{*}, b; z')$ is a $2$-connected rooted graph 
such that $V(B^{*}) \setminus \{s^{*}, b, z'\} \neq \emptyset$. 
Recall that 
$t_{0} \notin V(B) \setminus \{b, z'\}$ for the case where $\HF$ is of type~$3^{\fl}$. 
For a vertex $v \in V(B^{*}) \sm\{s^{*}, b, z'\}$, %
it follows from Claims~\ref{claim:at most l+1}, 
~\ref{claim:subsetSbx}~(1) 
and $\lp \ge 2$ that  
\[
\deg_{B^{*}}(v) \ge 
\begin{cases}
\deg_{G}(v) - |\{x\}| 
\ge 
k 
\ge 
(k-\lp+1) +1 
& \text{if $N_{G}(v) \cap S = \emptyset$,}\\
\deg_{G}(v) - \el + 1 
\ge 
\deg_{G}(v) - \lp + 1 
\ge
(k-\lp+1) +1 
& \text{otherwise,}
\end{cases} 
\]
and thus
$\delta(B^{*}, s^{*}, b ;z') 
\ge (k - \lp + 1) +1$.
By the induction hypothesis, 
$B^{*}$ contains $k - \lp +1$ admissible $(s^{*}, b)$-paths.
Therefore
$G[B \cup \big( N_{G}(B-b) \cap S\big)]$ 
contains $k-\lp+1$ admissible $(S, b)$-paths $\oraP{1}, \dots, \oraP{k-\lp+1}$. 
Let $s_{i} \in V(P_{i}) \cap S$ for $1 \le i \le k-\lp+1$, 
and let $\ora{R'}$ be a $(b, y)$-path internally disjoint from $B$ in $G - V(H^{\sh})$. 
Then 
$s_{i} \oraP{i} b \ora{R'} y$ is an $(s_{i}, y)$-path in $G[\big( V(G) \setminus V(H^{\sh}) \big) \cup \{s_{i}\}]$ 
for $1 \le i \le k-\lp+1$. 
On the other hand, 
it follows from 
Fact~\ref{fact:paths in H}~\ref{xS} 
and the definition of type~II 
that for each $s_{i}$, 
$H^{\sh}$ contains $\lp$ admissible $(x, s_{i})$-paths 
of lengths 
$2, \ldots, 2\lp$ (if $\HF$ is of type~I or type~III)
or
$2,3$ (if $\HF$ is of type~II). 
Hence by Fact~\ref{fact:Lemma3.2}, 
$G$ contains $k \ \big( = (k - \lp +1) + \lp - 1 \big)$ admissible $(x, y)$-paths, a contradiction. 
\qed

Since $S \not= \emptyset$,
it follows from Claim~\ref{claim:l=1}
that 
\begin{center}
$S^{\sh} = S$, that is, 
$\HF$ is of type~III,  \, $\lp=\el = |S| = 1$, say $S = \{s_{1}\}$,  \, $H^{\sh} = \HF$.
\end{center}
Then the following hold (note that $t_{0} \notin \NCT$,
since $t_{0}$ is a cut-vertex of $\CF$): 
\begin{gather}
\label{neighbor of B-b}
\text{$N_{G}(B - b) \cap V(\HF) = \{x, s_{1}\}$ \, (by Claim~\ref{claim:subsetSbx}~(2))}, \\
\label{le1}
\text{$|N_{G}(v) \cap V(\HF)| \le 1$ for each $v \in \NCT$ \, (by Claim~\ref{claim:at most l+1}})
\end{gather}

\begin{claim}
\label{claim:kge3 and EedgeTC}
(1) $k \ge 3$, 
and 
(2) if $z \in V(\CF)$, then $N_{G}(\TF) \cap V(\CF - y) \not= \emptyset$. 
\end{claim}

\proof 
By Claim~\ref{claim:S=1NCT}, 
we have $N_{G}(\TF) \cap V(\CF) \not=\emptyset$. 
Therefore, 
it follows from 
Fact~\ref{fact:paths in H}~\ref{xT} that $k \ge 3$. 
Thus (1) holds. 
To show (2), 
suppose that 
$z \in V(\CF)$ and 
$N_{G}(\TF) \cap V(\CF - y) = \emptyset$. 
Since $N_{G}(\TF) \cap V(\CF) \not=\emptyset$, 
we have 
$N_{G}(\TF) \cap V(\CF) = \{y\}$. 
Let $G' = G - V(\CF - y)$. 
Since 
$G$ and $\HF$ are $2$-connected, 
$z \in V(\CF - y)$ 
and 
$|V(\HF)| \ge 4$, 
it follows that 
$(G', x, y; s_{1})$ is a $2$-connected rooted graph 
such that  
$\delta(G', x, y; s_{1}) \ge \delta(G, x, y; z)$. 
Therefore, 
by the induction hypothesis, 
$G'$ (and also $G$) contains $k$ admissible $(x, y)$-paths in $G$, a contradiction. 
Thus (2) also holds. 
\qed

\begin{claim}
\label{claim:notfeasible}
(1) $V(B-b)\cap \{y,z'\}\not=\emptyset$,
and 
(2) $\CF$ is not  a block.
\end{claim}

\proof
Suppose that $y, z' \not\in V(B-b)$. 
(Note that then $B$ satisfies \ref{Cis2-conn} or \ref{existsFEB}-\ref{ynotinB}.) 
Recall that (\ref{neighbor of B-b}) holds. 
If 
$|N_{G}(s_{1}) \cap V(B)| \ge 2$, 
then let 
$B' = G[B \cup \{x, s_{1}\}]$ 
and 
$z_{B} = s_{1}$; 
if 
$|N_{G}(s_{1}) \cap V(B)| = 1$, 
say $N_{G}(s_{1}) \cap V(B) = \{v\}$, 
then let 
$B' = G[B \cup x]$ 
and 
$z_{B} = v$. 
Note that $|V(B)| \ge 3$, 
since 
$V(B) \setminus \{b, z'\} \neq \emptyset$ 
and $\delta(G, x, y; z) \ge k + 1 \ge 4$ by Claim~\ref{claim:kge3 and EedgeTC}~(1). 
Then 
$(B', x, b; z_{B})$ is a $2$-connected rooted graph and 
$\delta(B', x, b; z_{B}) \ge \delta(G, x, y; z)$. 
Hence by the induction hypothesis, 
$B'$ contains $k$ admissible $(x, b)$-paths $\oraP{1}, \dots, \oraP{k}$. 
Then 
$x \oraP{i} b \ora{R} y$ ($1 \le i \le k$) 
are $k$ admissible $(x, y)$-paths in $G$, a contradiction. 
Thus (1) holds.

Suppose next that 
$\CF$ is a block. 
Then by \ref{Cis2-conn}, note that $B = \CF$, $b=y$ and $z' = z$. 
In particular,
Claim~\ref{claim:notfeasible}~(1) implies that $z = z' \in V(\CF)$. 
Then 
by Claim~\ref{claim:kge3 and EedgeTC}~(2),
$N_{G}(\TF) \cap V(B - b) = N_{G}(\TF) \cap V(\CF - y) \not=\emptyset$. 
But, this contradicts Claim~\ref{claim:subsetSbx}~(1). 
Thus (2) also holds. 
\qed

Recall that $(B, b, z')$ denotes any one of 
\ref{Cis2-conn}, \ref{existsFEB}-\ref{yinB}, \ref{existsFEB}-\ref{ynotinB} and \ref{existsFB}. 
By Claim~\ref{claim:notfeasible}, 
$\CF$ has exactly two end-blocks, 
and 
each end-block of $\CF$ contains exactly one of $z$ and $y$ as a non-cut-vertex of $\CF$ 
(otherwise, there is a feasible end-block of $C$ which satisfies no Claim~\ref{claim:notfeasible}~(1)). 
In particular, 
$(\CF, z, y)$ is a $2$-connected rooted graph.

Let 
$B_{1}, \dots, B_{t}$ ($t \ge 2$) be all the blocks of $\CF$ 
such that 
$V(B_{i}) \cap V(B_{i + 1}) \neq \emptyset$ for $1 \le i \le t-1$, 
say $V(B_{i}) \cap V(B_{i + 1}) = \{b_{i}\}$ for $1 \le i \le t-1$. 
Without loss of generality, 
we may assume that 
$z \in V(B_{1}) \setminus \{b_{1}\}$ 
and 
$y \in V(B_{t}) \setminus \{b_{t-1}\}$, 
and 
let 
$b_{0} = z$ 
and 
$b_{t} = y$. 
Then $B = B_{p}$ for some $p$ with $1 \le p \le t$. 
Note that 
$B_{p}, b_{p-1}, b_{p}$ satisfy
\ref{existsFEB}-\ref{yinB} (if $p = t$) 
or
\ref{existsFEB}-\ref{ynotinB} (if $1=p<t$) 
or 
\ref{existsFB} (if $2 \le p < t$)
as $(B, b, z') = (B_{p}, b_{p}, b_{p-1})$, 
and so 
it follows from Claim~\ref{claim:subsetSbx}~(1) 
that 
\begin{align}
\label{neighbor of T}
N_{G}(\TF) \cap V(B_{p} - b_{p}) = \emptyset. 
\end{align}

\begin{claim}\label{claim:if Bp has order>=3}
$N_{G}(\TF) \cap \big(\bigcup_{1 \le i \le p-1}V(B_{i}) \big)= \emptyset$. 
%
\end{claim}

\proof 
Suppose that 
$N_{G}(\TF) \cap \big( \bigcup_{1 \le i \le p-1}V(B_{i}) \big) \neq \emptyset$. 
Then 
it follows from Fact~\ref{fact:paths in H}~\ref{xT} that 
$G\left[\HF  \cup \big(\bigcup_{1 \le i \le p-1}V(B_{i})\big) \right]$ 
contains two admissible $(x, b_{p-1})$-paths. 
On the other hand, 
since 
$v \in \NCT$ for $v \in V(B_{p}) \setminus \{b_{p-1}, b_{p}\}$, 
it follows from (\ref{le1}) that 
$(B_{p}, b_{p-1}, b_{p}; z)$ is a $2$-connected rooted graph 
such that 
$\delta(B_{p}, b_{p-1}, b_{p}; z) \ge \delta(G, x, y; z) - 1 \ge (k-1) + 1$, 
and hence 
the induction hypothesis 
yields that $B_{p}$ contains $k-1$ admissible $(b_{p-1}, b_{p})$-paths. 
Let 
$\ora{R'}$ be a $(b_{p}, y)$-path 
in $G[\bigcup_{p \le i \le t}V(B_{i})]$. 
Then 
$b_{p-1} \oraP{i} b_{p} \ora{R'} y$ ($1 \le i \le k-1$) are $k-1$ admissible $(b_{p-1}, y)$-paths. 
Therefore, by Fact~\ref{fact:Lemma3.2}, 
$G$ contains $k \ \big( = (k-1) + 2 - 1 \big)$ admissible $(x, y)$-paths, a contradiction. 
\qed

Choose $B = B_{p}$ so that $p$ is as large as possible. 
If $p=t$, 
then by (\ref{neighbor of T}) and Claim~\ref{claim:if Bp has order>=3}, 
we have 
$N_{G}(\TF) \cap V(\CF - y) = \emptyset$; 
since 
$z \in V(B_{1}) \setminus \{b_{1}\}$, 
this contradicts Claim~\ref{claim:kge3 and EedgeTC}~(2). 
Thus $p < t$ and the choice of $B = B_{p}$ implies that 
$|V(B_{t})| = 2$, i.e., $V(B_{t}) = \{b_{t-1}, y\} \ (= \{b_{t-1}, b_{t}\})$.

Recall that any core with respect to $(y, x)$ is of type $3$ 
(by \ref{refno}, \ref{renfomin}, Remark~\ref{remark:lcore}, Claim~\ref{claim:2-connected and xy notin E}). 
By (\ref{neighbor of B-b}), 
$\deg_{G}(x) \ge |\TF| + |N_{G}(x) \cap (B_{p} - b_{p})| \ge |\TF| + 1$, 
and so 
\ref{degxy} 
yields 
that 
$\deg_{G}(y) \ge |\TF| + 1$. 
Since 
$N_{G}(y) \subseteq \HF \cup b_{t-1}$, 
we obtain 
$|N_{G}(y) \cap V(\HF)| \ge |\TF| \ge 2$. 
This implies that 
$N_{G}(y) \cap V(\HF) = \TF$ 
and 
$N_{G}(b_{t-1}) \cap \TF=\emptyset$ 
(otherwise, $xy \in E(G)$ or there exists a core of type~1 with respect to $(y, x)$, 
a contradiction). 
If $p=t-1$, then 
by the same argument as the case $p=t$, 
we get a contradiction to Claim~\ref{claim:kge3 and EedgeTC}~(2). 
Thus $p < t - 1$ and 
the choice of $B = B_{p}$ implies that 
$|V(B_{t-1})| = 2$, i.e., $V(B_{t-1}) = \{b_{t-2}, b_{t-1}\}$. 
Since $\deg_{G}(y) = |\TF| + 1$,
we have $N_{G}(x) = \TF \cup \big( N_{G}(x) \cap (B_{p} - b_{p}) \big)$, 
and so $x \notin N_{G}(b_{t-1})$ because $b_{t-1} \notin V(B_{p})$. 
Therefore $N_{G}(b_{t-1})\subseteq \{y,b_{t-2},s_{1}\}$.
Since $\deg_{G}(b_{t-1}) \ge k+1$,
we obtain $k \le 2$,
contradicting to Claim~\ref{claim:kge3 and EedgeTC}~(1).

This completes the proof of Theorem~\ref{thm:three vertices in rooted graph}. 
\qed

\medskip

We finally give the proof of Theorem~\ref{thm:generalization of BV}.

\medskip
\noindent
\textit{Proof of Theorem~\ref{thm:generalization of BV}.}~It suffices to show the case where 
a given graph is connected. 
Let $k \ge 2$ be an integer, 
and let $G$ be a connected graph of order at least three 
having at most two vertices of degree less than $k + 1$. 
Let 
$x$ and $z$ be two vertices of degree less than $k + 1$ if exist; 
otherwise, 
let $x$ and $z$ be arbitrary two vertices. 
Suppose now that $G$ is a counterexample.

We first consider the case where $G$ is $2$-connected. 
Choose arbitrary edge $xy$ in $G$ 
(possibly $y = z$). 
Since $|V(G)| \ge 3$ and $\deg_{G}(v) \ge k + 1 \ge 3$ for $v \in V(G) \setminus \{x, z\}$, 
we have $V(G) \setminus \{x, y, z\} \neq \emptyset$ 
and 
$\deg_{G}(v) \ge k + 1$ for $v \in V(G) \setminus \{x, y, z\}$. 
Hence 
Theorem~\ref{thm:three vertices} yields that 
$G$ contains $k$ admissible $(x, y)$-paths. 
By adding $xy$ to each of the $k$ paths, 
we obtain $k$ admissible cycles, a contradiction. 
Thus $G$ is not $2$-connected.

Suppose that 
there exists an end-block $B$ with cut-vertex $b$ 
such that 
$|V(B)| \ge 3$ and $|V(B - b) \cap \{x, z\}| \le 1$. 
Let $x' \in V(B - b) \cap \{x, z\}$ if exists; 
otherwise, $x' \in V(B-b)$. 
Then the same argument as in the case where $G$ is $2$-connected 
can work with $(G, x, z) = (B, x', b)$, 
and so we obtain $k$ admissible cycles in $B$, a contradiction. 
This, together with the degree condition, implies that 
the block-tree of $G$ is a path, 
and 
the two end-blocks of $G$ are the $x$-end-block and the $z$-end-block, respectively. 
Since 
$|V(G)| \ge 3$ and $\deg_{G}(v) \ge k + 1 \ge 3$ for $v \in V(G) \setminus \{x, z\}$, 
there exists a block $B$ with exactly two cut-vertices $b_{1}, b_{2}$ 
such that $|V(B)| \ge 3$. 
Then by replacing 
$(G, x, z)$ and $(B, b_{1}, b_{2})$ in the above argument for the case where $G$ is $2$-connected, 
we obtain $k$ admissible cycles in $B$, a contradiction again. 
\qed



\end{document}